\documentclass[twocolumn]{autart}

\usepackage{graphicx}
\usepackage{epsfig} 
\usepackage{times} 
\usepackage{amsmath}
\usepackage{amssymb}  
\usepackage{color}
\usepackage{amsfonts}
\usepackage{subfigure}
\usepackage{multirow}
\usepackage{multicol}
\usepackage{siunitx}
\usepackage{overpic}
\usepackage{mathrsfs}  
\usepackage{wrapfig}

\usepackage{booktabs}
\usepackage{threeparttable}

\usepackage{epstopdf}
\newtheorem{example}{Example}
\newtheorem{lemma}{Lemma}

\newtheorem{theorem}{Theorem}

\newtheorem{assumption}{Assumption}

\usepackage{xspace}

\begin{document}


\begin{frontmatter}



\title{Compensation of Actuator Dynamics Governed by \\ Quasilinear Hyperbolic PDEs}

\author[nik]{Nikolaos Bekiaris-Liberis\thanksref{footnoteinfo}}\ead{nikos.bekiaris@dssl.tuc.gr} and  
\author[miros]{Miroslav Krstic}\ead{krstic@ucsd.edu}

\address[nik]{Department of Production Engineering \& Management, Technical University of Crete, Chania, 73100, Greece}

\address[miros]{Department of Mechanical \& Aerospace Engineering, University of California, San Diego, La Jolla, CA 92093-0411, USA}            

\thanks[footnoteinfo]{Corresponding author.}



\maketitle
\thispagestyle{empty}     
\pagestyle{empty}

\begin{abstract}
We present a methodology for stabilization of general nonlinear systems with actuator dynamics governed by general, quasilinear, first-order hyperbolic PDEs. Since for such PDE-ODE cascades the speed of propagation depends on the PDE state itself (which implies that the prediction horizon cannot be a priori known analytically), the key design challenge is the determination of the predictor state. We resolve this challenge and introduce a PDE predictor-feedback control law that compensates the transport actuator dynamics. Due to the potential formation of shock waves in the solutions of quasilinear, first-order hyperbolic PDEs (which is related to the fundamental restriction for systems with time-varying delays that the delay rate is bounded by unity), we limit ourselves to a certain feasibility region around the origin and we show that the PDE predictor-feedback law achieves asymptotic stability of the closed-loop system, providing an estimate of its region of attraction. Our analysis combines Lyapunov-like arguments and ISS estimates. Since it may be intriguing as to what is the exact relation of the cascade to a system with input delay, we highlight the fact that the considered PDE-ODE cascade gives rise to a system with input delay, with a delay that depends on past input values (defined implicitly via a nonlinear equation). 


\end{abstract}

\end{frontmatter}

\interdisplaylinepenalty=2500 



\section{Introduction}
\subsection{Motivation}
Numerous processes may be described by quasilinear, first-order hyperbolic Partial Differential Equations (PDEs) cascaded with nonlinear Ordinary Differential Equations (ODEs), such as, for example, communication networks \cite{espita}, blood flow \cite{borch1}, sewer networks \cite{lunge}, production systems \cite{gotligh}, vehicular traffic flow \cite{herty}, piston dynamics \cite{latanzo}, and automotive engines \cite{dep}, \cite{jankovic1}, \cite{jankovic2} to name only a few \cite{richard}. Despite their popularity, despite the fact that predictor-based control laws now exist for nonlinear systems with input delays that may depend on the ODE state \cite{bek4}, \cite{bek5}, \cite{bek6}, \cite{bek7}, \cite{cai1}, \cite{cai2} as well as the uncontrolled- or controlled-boundary value of the PDE state \cite{delph5}, \cite{delph1}, \cite{delph4}, \cite{diagne1}, and despite the existence of several results on boundary stabilization of quasilinear, first-order hyperbolic PDEs, such as, for example, \cite{blandin}, \cite{coron}, \cite{hu}, \cite{krstic burgers}, \cite{prier}, \cite{vazquez1}, \cite{vazquez2}, no result exists on the compensation of actuator dynamics governed by quasilinear, first-order hyperbolic PDEs for nonlinear systems.


\subsection{Contributions}
In this paper, we consider the problem of stabilization of nonlinear ODE systems through transport actuator dynamics governed by quasilinear, first-order hyperbolic PDEs. We develop a novel PDE predictor-feedback law, which compensates the PDE actuator dynamics. Since the speed of propagation depends on the PDE state itself, the key idea in our design is the construction of the PDE predictor state. This construction is by far non-trivial and cannot follow in a straightforward way employing the results from \cite{diagne1}, which is perhaps the only available work dealing with the problem of complete compensation of an input-dependent input delay (note that the designs in \cite{delph1}, \cite{delph4}, \cite{delph5}, don't aim at achieving complete delay compensation). The reason is that the transport speed in the class of systems considered in \cite{diagne1} depends only on the {\em uncontrolled-boundary} value of the PDE state rather than on the PDE state itself, as it is the case here. 

Furthermore, we show that the PDE predictor-feedback design achieves local asymptotic stability in the $C^1$ norm of the actuator state. The reason for obtaining only a regional result, restricting the $C^1$ norm of the PDE state, is the possibility of appearance of multivalued solutions, or, in other words, the appearance of shock waves, in the solutions of quasilinear, first-order hyperbolic PDEs. We show, within our stability analysis, that this issue is avoided, limiting the $C^1$ norm of the solutions and the initial conditions. This limitation may alternatively be expressed as the fundamental limitation in stabilization of systems with time-varying input delays that the delay rate is bounded by unity--for the class of systems considered here, giving rise to an input delay that depends on the actuator state and its derivative, the satisfaction of this restriction is guaranteed by confining the size of the actuator state and its derivative. The proof of asymptotic stability in the $C^1$ norm of the actuator state is established employing Lyapunov-like arguments as well as Input-to-State Stability (ISS) estimates. 

In order to make the presentation of our control design methodology accessible to both readers who are experts on PDEs and readers who are experts on delay systems we highlight the relation of the PDE-ODE cascade to a system with input delay that is defined implicitly through a nonlinear equation, which involves the input value at a time that depends on the delay itself, and, moreover, we present the predictor-feedback design in this representation as well. 

\subsection{Organization}
We start in Section \ref{Sec1} where we present the class of systems under consideration as well as the PDE predictor-feedback control design. We provide an alternative, delay system representation of the considered PDE-ODE cascade in Section \ref{Sec2}. In Section \ref{sec3} we prove the local asymptotic stability of the closed-loop system under the proposed controller. Concluding remarks are provided in Section \ref{sec con}.

{\em Notation:} We use the common definition of class $\mathcal{K}$, $\mathcal{K}_{\infty}$ and $\mathcal{KL}$ functions from \cite{khalil}. For an $n$-vector, the norm $|\cdot|$ denotes the usual Euclidean norm. For a scalar function $u\in C[0,1]$ we denote by $\|u(t)\|_{\infty}$ its respective maximum norm, i.e., $\|u(t)\|_{\infty}=\max_{x\in[0,1]}|u(x,t)|$. For a scalar function $u_x\in C[0,1]$ we denote by $\|u_x(t)\|_{\infty}$ its respective maximum norm, i.e., $\|u_x(t)\|_{\infty}=\max_{x\in[0,1]}|u_x(x,t)|$. For a vector valued function $p\in C[0,1]$ we denote by $\|p(t)\|_{\infty}$ its respective maximum norm, i.e., $\|p(t)\|_{\infty}=\max_{x\in[0,1]}\sqrt{p_1(x,t)^2+\ldots +p_n(x,t)^2}$. For a vector valued function ${p}_x\in C[0,1]$ we denote by $\|{p}_x(t)\|_{\infty}$ its respective maximum norm, i.e., $\|{p}_x(t)\|_{\infty}=\max_{x\in[0,1]}\sqrt{{{p}_1}_x(x,t)^2+\ldots+ {{p}_n}_x(x,t)^2}$. We denote by $C^j(A;E)$ the space of functions that take values in $E$ and have continuous derivatives of order $j$ on $A$.

\section{Problem Formulation and Predictor-Feedback Control Design}
\label{Sec1}
We consider the following system
\begin{eqnarray}
\dot{X}(t)&=&f\left(X(t),u(0,t)\right)\label{pdesys1}\\
u_t(x,t)&=&v\left(u(x,t)\right)u_x(x,t)\label{pdesys3}\\
u\left(1,t\right)&=&U(t),\label{pdesys2}
\end{eqnarray}
where $X\in\mathbb{R}^n$ and $u\in\mathbb{R}$ are ODE and PDE states, respectively, $t\geq0$ is time, $x\in[0,1]$ is spatial variable, $U$ is control input, and $f:\mathbb{R}^n\times\mathbb{R}\to\mathbb{R}^n$ is a continuously differentiable vector field that satisfies $f(0,0)=0$.

The following assumptions are imposed on system (\ref{pdesys1})--(\ref{pdesys2}). 
\begin{assumption}
\label{ass D}
\rm{Function $v:\mathbb{R}\to\mathbb{R}_+$ is twice continuously differentiable and there exists a positive constant $\underline{v}$ such that the following holds
\begin{eqnarray}
v\left(u\right)\geq\underline{v},\quad\mbox{for all $u\in\mathbb{R}$}.\label{assv}
\end{eqnarray}}
\end{assumption}
\begin{assumption}
\label{assf}
\rm{System $\dot{X}=f\left(X,\omega\right)$ is strongly forward complete with respect to $\omega$.}
\end{assumption}
\begin{assumption}
\label{assi}
\rm{There exists a twice continuously differentiable feedback law $\kappa:\mathbb{R}^n\to\mathbb{R}$, with $\kappa(0)=0$, which renders system $\dot{X}=f\left(X,\kappa(X)+\omega\right)$ input-to-state stable with respect to $\omega$.}
\end{assumption}

Assumption \ref{ass D} is a prerequisite  for the well-posedness of the predictor state, which is defined in the next paragraph. It guarantees that transport is happening only in the direction away from the input, or, in other words (see also the discussion in the next section), it ensures that the input delay is positive as well as uniformly bounded. Assumption \ref{assf} (see, e.g., \cite{angeli}) and Assumption \ref{assi} (see, e.g., \cite{sontag}) are standard ingredients of the predictor-feedback control design methodology (see, e.g., \cite{bek5}, \cite{krstic book}, \cite{krstic}). The former implies that the state $X$ of system (\ref{pdesys1}) doesn't escape to infinity before the control signal $U$ reaches it, no matter the size of the delay (see, e.g., \cite{bek5}, \cite{krstic book}, \cite{krstic}), while the latter guarantees the existence of a nominal feedback law that renders system (\ref{pdesys1}) input-to-state stable in the absence of the transport actuator dynamics (i.e., in the absence of the input delay).

The predictor-feedback control law for system (\ref{pdesys1})--(\ref{pdesys2}) is given by
 \begin{eqnarray}
U(t)=\kappa\left(p\left(1,t\right)\right),\label{controller}
\end{eqnarray}
where for all $t\geq0$
\begin{eqnarray}
p\left(x,t\right)&=&X(t)+\int_0^{x}f\left(p(y,t),u(y,t)\right)\nonumber\\
&&\times\Gamma\left(u(y,t),u_y(y,t),y\right)dy,\quad x\in[0,1]\label{p1}
\end{eqnarray}
with\footnote{Note that $\Gamma$ can be written as $\Gamma\left(u(x,t),u_x(x,t),x\right)=\frac{\partial \frac{x}{v\left(u(x,t)\right)}}{\partial x}$.}
\begin{eqnarray}
\Gamma\left(u(x,t),u_x(x,t),x\right)&=&\frac{1}{v\left(u(x,t)\right)}\nonumber\\
&&-\frac{xv'\left(u(x,t)\right)u_x(x,t)}{v\left(u(x,t)\right)^2},\nonumber\\
&& x\in[0,1].\label{gamma}
\end{eqnarray}
For implementing the predictor-feedback law (\ref{controller})--(\ref{gamma}), besides measurements of the ODE state $X(t)$ and the PDE state $u(x,t)$, $x\in[0,1]$, for all $t\geq0$, the availability of the spatial derivative of $u$, namely, $u_x(x,t)$, $x\in[0,1]$, for all $t\geq0$, is required. The latter may be obtained either via direct measurements of $u_x$ or by a numerical computation of $u_x$, employing the measurements of $u$. The implementation and approximation problems of predictor-feedback control laws are tackled, for example, in \cite{kar1}, \cite{mondie}, \cite{zhong}.

In order to guarantee the well-posedness of the predictor state (\ref{p1}) and the system the following feasibility condition on the closed-loop solutions and the initial conditions needs to be satisfied 
\begin{eqnarray}
-M<\frac{v'\left(u(x,t)\right)u_x(x,t)}{v\left(u(x,t)\right)}<1,\nonumber\\
 \mbox{for all $x\in[0,1]$ and $t\geq0$},\label{condition strictn1}
\end{eqnarray}
for some $M>0$. In the next section we provide some explanatory remarks on the feasibility condition (\ref{condition strictn1}) and Assumption \ref{ass D}, capitalizing on the relation of the PDE-ODE cascade (\ref{pdesys1})--(\ref{pdesys2}) to a system with a delayed-input-dependent input delay.

\begin{example}
\rm{To illustrate the control design and its implementation we present here a rather pedagogical example, which results in a predictor-feedback law defined explicitly in terms of $X$, $u$, and $u_x$. Consider an unstable, scalar linear system with actuator dynamics governed by a quasilinear, first-order hyperbolic PDE given by
\begin{eqnarray}
\dot{X}(t)&=&X(t)+u(0,t)\label{exsys1}\\
u_t(x,t)&=&\left(u(x,t)^2+1\right)u_x(x,t)\\
u(1,t)&=&U(t).\label{exsys3}
\end{eqnarray}
System (\ref{exsys1})--(\ref{exsys3}) satisfies all of the Assumptions \ref{ass D}--\ref{assi}  and a nominal control law may be chosen as $U(t)=-2X(t)$. Thus, the predictor-feedback control law is given by
\begin{eqnarray}
U(t)=-2p(1,t),\label{conex1}
\end{eqnarray}
where, exploiting the fact that $\Gamma=\frac{\partial \frac{x}{u(x,t)^2+1}}{\partial x}$ as well as the linearity of the system, the predictor state $p$, defined in (\ref{p1}), may be written in the present case as\footnote{To see this note that, for the case of system (\ref{exsys1})--(\ref{exsys3}), the predictor state $p$ in (\ref{p1}) satisfies, for each $t$, the ODE in $x$ given as $p_x(x,t)=\left(p(x,t)+u(x,t)\right)\frac{\partial \frac{x}{u(x,t)^2+1}}{\partial x}$, with initial condition $p(0,t)=X(t)$. Thus, solving this initial-value problem with respect to $x$ we obtain $p(x,t)=e^{\int_0^x\frac{\partial \frac{y}{u(y,t)^2+1}}{\partial y}dy}X(t)+\int_0^xe^{\int_y^x\frac{\partial \frac{r}{u(r,t)^2+1}}{\partial r}dr}u(y,t)\frac{\partial \frac{y}{u(y,t)^2+1}}{\partial y}dy$. Expression (\ref{p example}) then follows evaluating the integral in the first term of this relation and employing one step of integration by parts in the integral in the second.}
\begin{eqnarray}
p(x,t)&=&e^{\frac{x}{u(x,t)^2+1}}\left(X(t)+u(0,t)\vphantom{\int_0^x e^{-\frac{y}{u(y,t)^2+1}}u_y(y,t)dy}\right.\nonumber\\
&&\left.+\int_0^x e^{-\frac{y}{u(y,t)^2+1}}u_y(y,t)dy\right)-u(x,t),\nonumber\\
&& x\in[0,1].\label{p example}
\end{eqnarray}
For the numerical computation of the integral in (\ref{p example}) we employ a simple composite left-endpoint rectangular rule, where the spatial derivates of $u$ are numerically computed utilizing a forward finite difference scheme. We choose the initial conditions as
\begin{eqnarray}
u(x,0)&=&1,\quad\mbox{for all $x\in[0,1]$}\label{initial ex1}\\
X(0)&=& -0.7.\label{initial ex2}
\end{eqnarray}
 In Fig. \ref{fig1} we show the response of the system, whereas in Fig. \ref{fig2} we show the control effort. 
\begin{figure}[t]
\centering
\includegraphics[width=\linewidth]{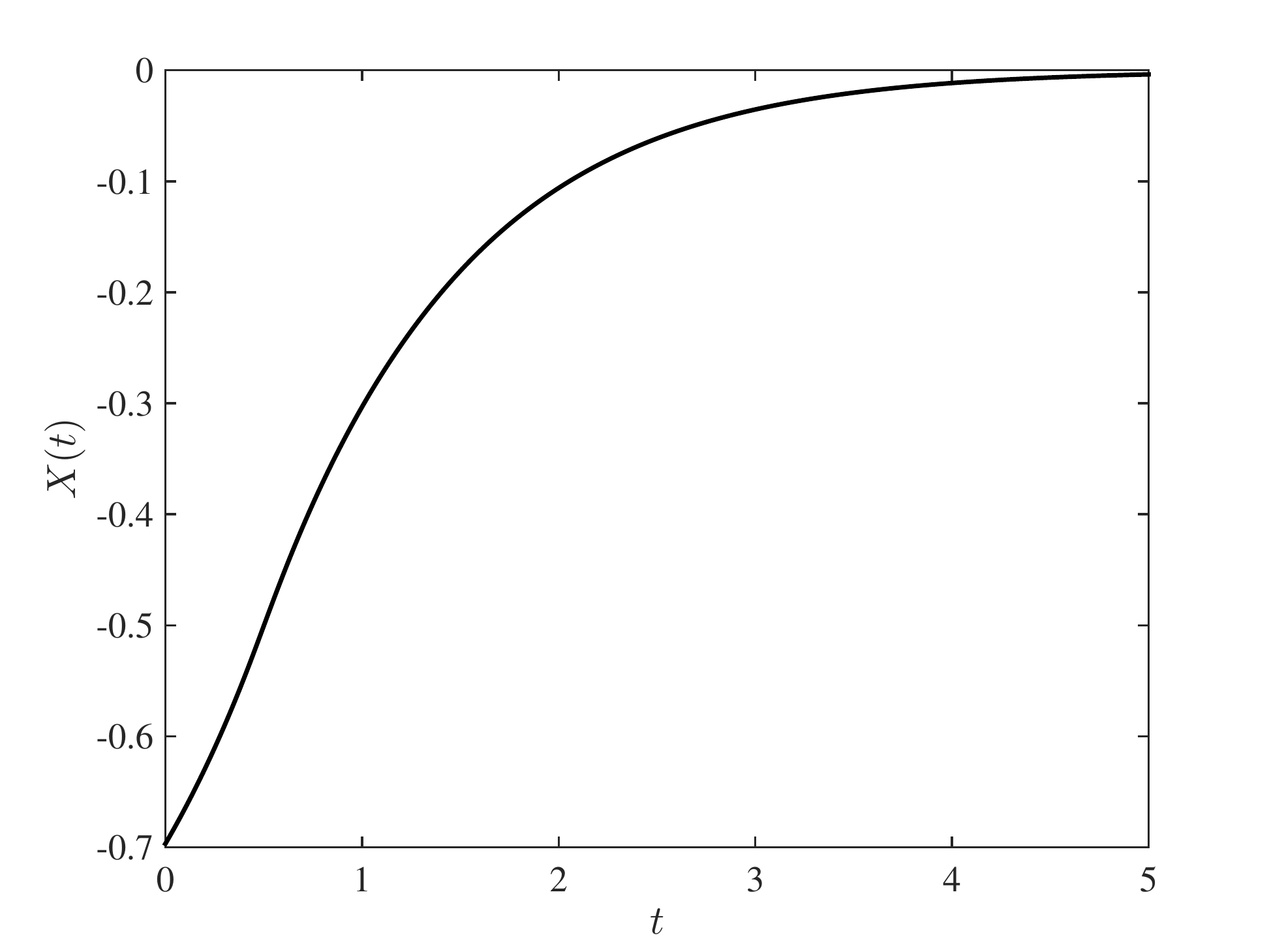}
\includegraphics[width=\linewidth]{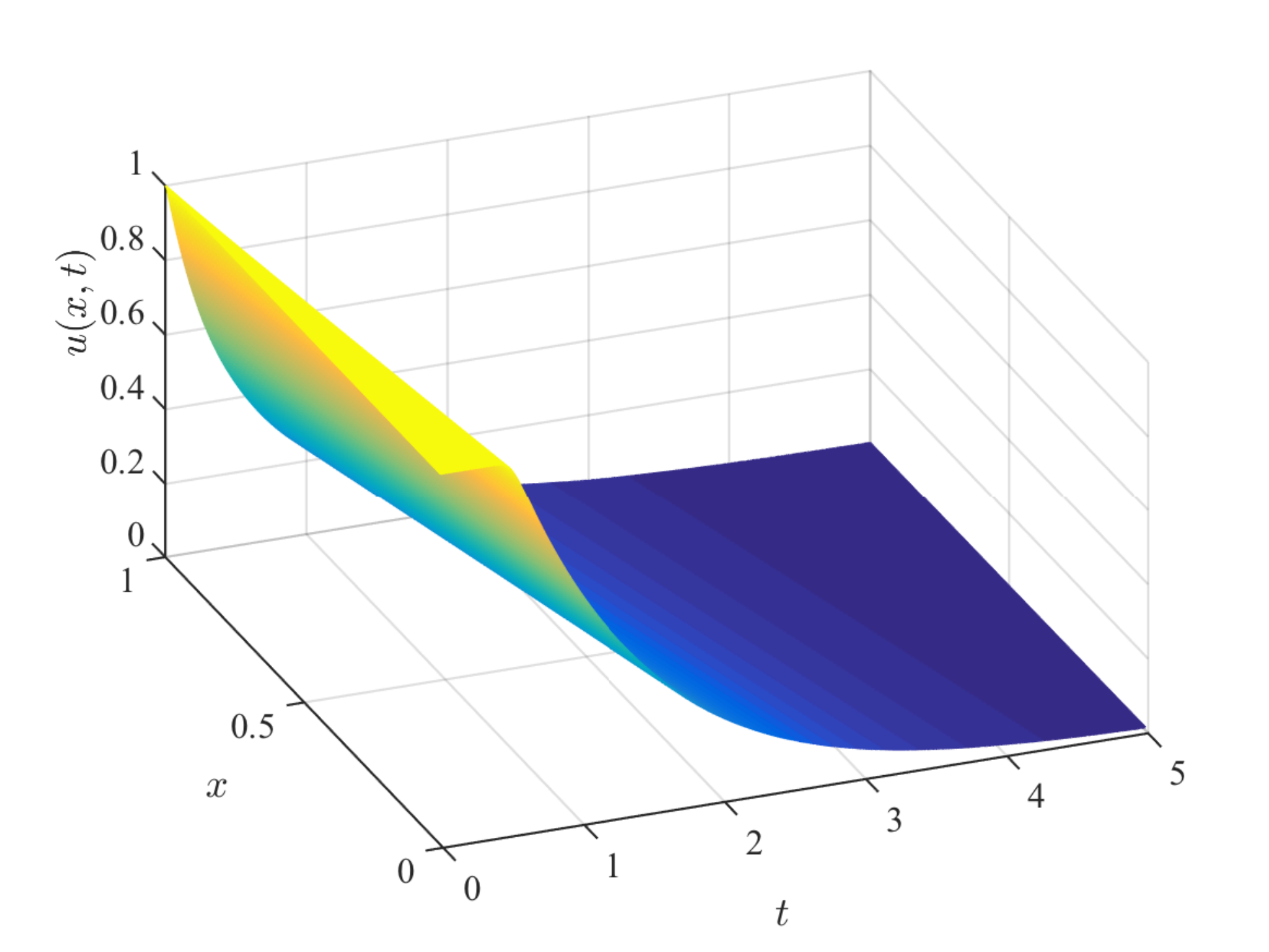}
\caption{Response of system (\ref{exsys1})--(\ref{exsys3}) with initial conditions (\ref{initial ex1}), (\ref{initial ex2}) under the predictor-feedback law (\ref{conex1}), (\ref{p example}).}
\label{fig1}
\end{figure}
\begin{figure}[t]
\centering
\includegraphics[width=\linewidth]{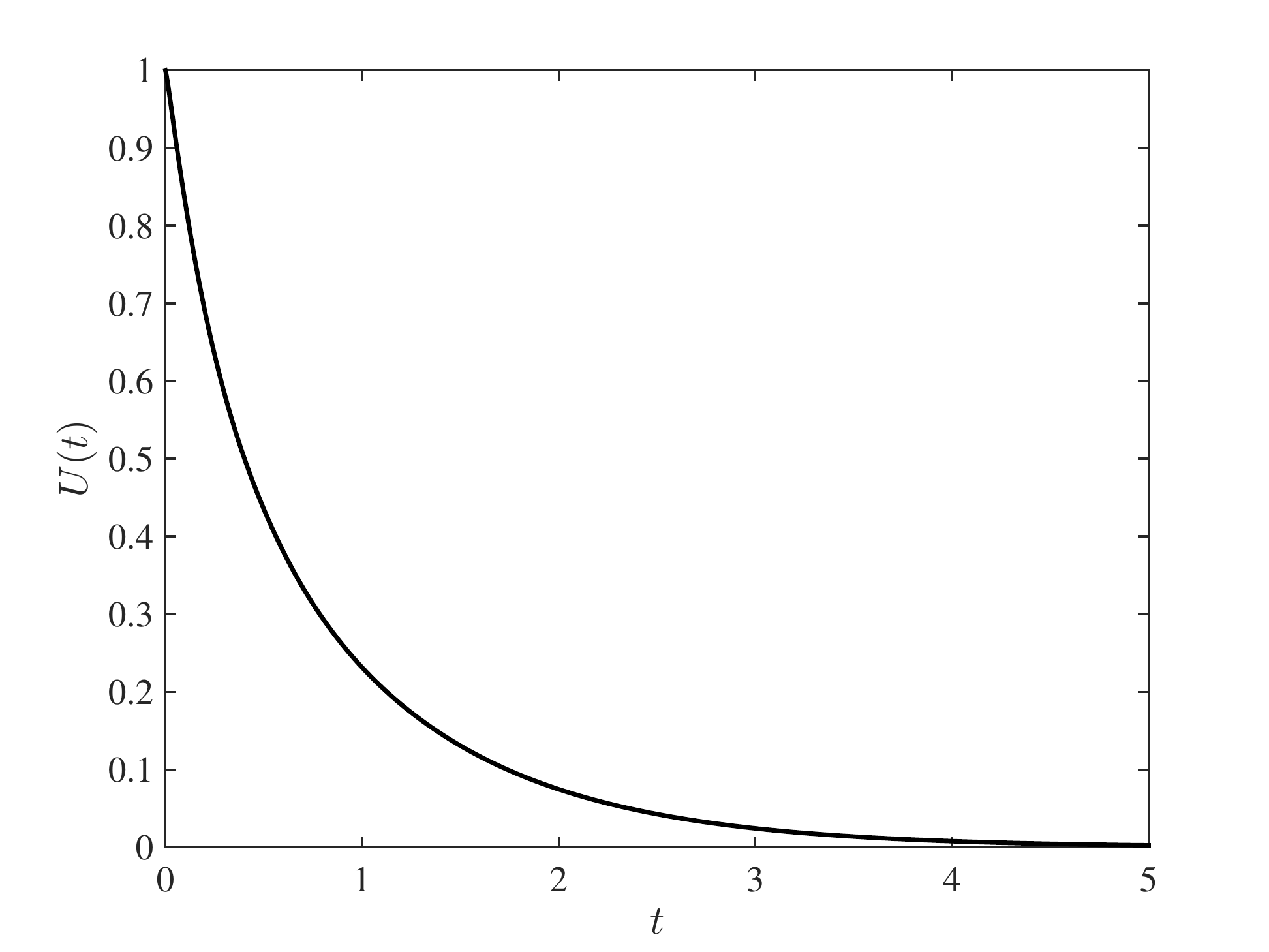}
\caption{Control effort (\ref{conex1}), (\ref{p example}).}
\label{fig2}
\end{figure}}
\end{example}

\section{Relation to a System with Delayed-Input-Dependent Input Delay}
\label{Sec2}
In this section, we highlight the fact that the PDE-ODE cascade (\ref{pdesys1})--(\ref{pdesys2}) may be viewed as a nonlinear system with an input delay. The fact that the transport speed depends on the PDE state itself, gives rise to a delay that is defined implicitly through a nonlinear equation, which incorporates the value of the input at a time that depends on the delay itself. 

The reasons for emphasizing this alternative representation of system (\ref{pdesys1})--(\ref{pdesys2}) are not merely pedagogical. Capitalizing on this relation, enables both, readers who are experts on PDEs and readers who are experts on delay systems, to digest the key conceptual ideas as well as the technical intricacies of our design and analysis methodologies, such as, for example, to better understand some of the inherent limitations of the stabilization problem for such systems (see Section \ref{subb1}). Moreover, this alternative point of view, offers to the designer two alternative control law representations (see Section \ref{subb2}), which may be very useful since, depending on the specific application, one representation may be more descriptive of the actual physical process as well as more suitable for implementation than the other (consider, for example, the case of control of traffic flow versus the case of control over a network).


\subsection{Derivation of the Delayed and Prediction Times}
Employing the method of characteristics (for details, see, e.g., \cite{courant}), it can be shown, see, e.g., \cite{petit1}, that the following holds
\begin{eqnarray}
u(0,t)=U\left(t-\frac{1}{v\left(u(0,t)\right)}\right).
\end{eqnarray}
Thus, defining the delayed time $\phi$, i.e., the time at which the value of the control signal $U$ that currently affects the system, namely, $u(0,t)$, was actually applied, as  
\begin{eqnarray}
\phi(t)=t-\frac{1}{v\left(u(0,t)\right)},\label{delay0}
\end{eqnarray}
we re-write system (\ref{pdesys1})--(\ref{pdesys2}) as 
\begin{eqnarray}
\dot{X}(t)=f\left(X(t),U\left(\phi(t)\right)\right),\label{plant}
\end{eqnarray}
where $\phi$ is defined implicitly, for all $t\geq0$, through relation
\begin{eqnarray}
\phi(t)=t-\frac{1}{v\left(U\left(\phi(t)\right)\right)}.\label{delay}
\end{eqnarray}
The prediction time $\sigma$, i.e., the time at which the value of the control signal $U$ currently applied, namely, $U(t)=u(D,t)$, will actually reach the system, is defined as the inverse function of $\phi$, namely,
 \begin{eqnarray}
 \sigma(t)&=&t+\frac{1}{v\left(U(t)\right)}.\label{sigma}
\end{eqnarray}
The invertibility of $\phi$ is guaranteed when the derivative of (\ref{delay}), given by
\begin{eqnarray}
\dot{\phi}(t)=\frac{1}{1-\frac{v'\left(U\left(\phi(t)\right)\right)U'\left(\phi(t)\right)}{v\left(U\left(\phi(t)\right)\right)^2}},\quad\mbox{for all $t\geq0$},
\end{eqnarray}
is positive for all times, or, equivalently, when the derivative of (\ref{sigma}), given by 
\begin{eqnarray}
\dot{\sigma}\left(t\right)={1-\frac{v'\left(U\left(t\right)\right)}{v\left(U\left(t\right)\right)^2}U'\left(t\right)},\quad\mbox{for all $t\geq0$}.\label{phi sigma}
\end{eqnarray}
is positive for all times.

\subsection{Interpretation of Assumption \ref{ass D} and Condition (\ref{condition strictn1})}
\label{subb1}
From (\ref{delay}) it is evident that the positivity assumption of $v$ guarantees that the delay is always positive, i.e., it guarantees the causality of system (\ref{plant}), and thus, also of system (\ref{pdesys1})--(\ref{pdesys2}). Moreover, relation (\ref{assv}) guarantees the boundness of the delay, i.e., it guarantees that the control signal eventually reaches the plant (\ref{plant}), and thus, also (\ref{pdesys1}). 

The interpretation of condition (\ref{condition strictn1}) is less obvious. When the derivative of the prediction (or the delayed) time is bounded and strictly positive both the prediction and delayed times are well-defined. Via (\ref{pdesys2}), it is evident from (\ref{phi sigma}) that this requirement is satisfied when condition (\ref{condition strictn1}) holds. In fact, condition (\ref{condition strictn1}) guarantees that the quasilinear first-order hyperbolic PDE (\ref{pdesys3}), (\ref{pdesys2}) exhibits smooth solutions and that the appearance of shock waves is avoided.  

To see this, note that when the right-hand side of (\ref{condition strictn1}) is violated the derivative of the delayed time becomes infinite (or, equivalently, the derivative of the prediction time becomes zero), that is, the delay disappears instantaneously (with slope approaching negative infinity). This implies that the delayed time becomes a multivalued function, which in turn is related to loss of regularity of the solutions to (\ref{pdesys3}), (\ref{pdesys2}) and the formation of a shock wave. From that point and on, the delayed time becomes a decreasing function, and thus, the plant receives all the more older information than already received (despite the fact that the direction of transport remains leftward since the transport speed is always positive).

Moreover, when $u$ is bounded, the regularity assumption on $v$ implies that the left-hand side of (\ref{condition strictn1}) may be violated when $u_x$ reaches negative infinity. In terms of the delay representation, it guarantees that the time derivative of the prediction time cannot become infinite, and thus, the predictor state remains well-posed.

\subsection{Predictor-Feedback Control Design for the Equivalent Delay System}
\label{subb2}
Defining 
\begin{eqnarray}
F\left(U\right)=\frac{1}{v\left(U\right)},
\end{eqnarray}
the predictor-feedback control law for system (\ref{plant}) with an input delay defined via (\ref{delay}) is given by
\begin{eqnarray}
U(t)=\kappa\left(P(t)\right),\label{con}
\end{eqnarray}
where the predictor $P$ is given for all $t\geq0$ by 
\begin{eqnarray}
P(\theta)&=&X(t)+\int_{\phi(t)}^{\theta}\left(1+F'\left(U(s)\right)\dot{U}(s)\right)\nonumber\\
&&\times f\left(P(s),U(s)\right)ds,\quad \mbox{for all $\phi(t)\leq\theta\leq t$}.\label{predictor}
\end{eqnarray}

The predictor-feedback control law (\ref{predictor}) is implementable since, for all $t\geq0$, it depends on the history of $U(s)$, over the window $\phi(t)\leq s\leq t$, the ODE state $X(t)$, which are assumed to be measured for all $t\geq0$, as well as on $\dot{U}(s)$, over the window $\phi(t)\leq s\leq t$, which is assumed to either be measured directly or computed from the values of $U(s)$, $\phi(t)\leq s\leq t$. Moreover, the implementation of the predictor-feedback design requires the computation at each time step of the delayed time $\phi$. This can either be performed by numerically solving relation (\ref{delay}), using the history of the actuator state, or by employing the following integral equation 
\begin{eqnarray}
\phi\left(\theta\right)&=&t-\int_{\theta}^{\sigma(t)}\frac{ds}{1+F'\left(U\left(\phi(s)\right)\right)U'\left(\phi(s)\right)},\nonumber\\
&&\mbox{for all $t\leq\theta\leq \sigma(t)$},\label{phi}
\end{eqnarray}
where $\sigma$ is defined in (\ref{sigma}).
The issue of implementation and approximation of nonlinear predictor feedbacks is addressed in detail in \cite{kar alone}, \cite{karr}, \cite{kar1}.

\section{Stability Analysis}
\label{sec3}
\begin{theorem}
\label{thm1}
Consider the closed-loop system consisting of the plant (\ref{pdesys1})--(\ref{pdesys2}) and the control law (\ref{controller})--(\ref{gamma}). Under Assumptions \ref{ass D}, \ref{assf}, and \ref{assi}, there exist a positive constant $\delta$ and a class $\mathcal{KL}$ function $\beta$ such that for all initial conditions $X(0)\in\mathbb{R}^n$ and $u(\cdot,0)\in C^1\left[0,1\right]$ which satisfy
\begin{eqnarray}
|X(0)|+\|u(0)\|_{\infty}+\|u_x(0)\|_{\infty}<\delta,\label{region}
\end{eqnarray}
as well as the compatibility conditions 
\begin{eqnarray}
u\left(1,0\right)&=&\kappa\left(p\left(1,0\right)\right)\label{comb1}\\
u_x\left(1,0\right)&=&\frac{\partial \kappa\left(p\left(1,0\right)\right)}{\partial p} f\left(p\left(1,0\right),u\left(1,0\right)\right)\nonumber\\
&&\times\Gamma\left(u(1,0),u_x(1,0),1\right)\label{comb2},
\end{eqnarray}
there exists a unique solution to the closed-loop system with $X(t)\in C^1[0,\infty)$, $u(x,t)\in C^1\left([0,1]\times [0,\infty)\right)$, and the following holds
\begin{eqnarray}
\Omega(t)&\leq&\beta\left(\Omega(0),t\right),\quad\mbox{for all $t\geq0$}\label{estimate}\\
\Omega(t)&=&|X(t)|+\|u(t)\|_{\infty}+\|u_x(t)\|_{\infty}.\label{omeg}
\end{eqnarray}
\end{theorem}

The proof of Theorem \ref{thm1} is based on the following lemmas whose proofs can be found in Appendix A. 

\begin{lemma}
\label{lemma2}
The variable
\begin{eqnarray}
u(x,t)-\kappa\left(p(x,t)\right)=w(x,t),\label{w1}
\end{eqnarray}
where $p$ is defined in (\ref{p1}), satisfies
\begin{eqnarray}
w_t(x,t)&=&v\left(u(x,t)\right)w_x(x,t)\label{pdent}\\
w(1,t)&=&0.\label{pden1t}
\end{eqnarray}
Moreover, system (\ref{pdesys1}) can be written as
\begin{eqnarray}
\dot{X}(t)&=&f\left(X(t),\kappa\left(X(t)\right)+w(0,t)\right)\label{pdesysn1t}.
\end{eqnarray}
\end{lemma}


Note that, differently with previous work on predictor-feedback design, the variable $w$ is just viewed as a new variable, which is expressed in terms of the state $(X,u)$ via (\ref{w1}), (\ref{p1}), (\ref{gamma}), rather than as a transformation of the original state $u$. Thus, an inverse transformation is not required, which doesn't affect the analysis (see Lemma \ref{lemma5} below and its proof in Appendix A). The reason for this alternative point of view is that the expression for the potential inverse transformation would require the definition of an alternative, rather complex representation of the predictor state $p$ that would depend on the new variable $w$, which would add unnecessary complexity in the analysis.

The next lemma establishes an asymptotic stability estimate for state variables $\left(X,w(x)\right)$, $x\in[0,1]$, exploiting the cascade structure of system (\ref{pdent})--(\ref{pdesysn1t}).

\begin{lemma}
\label{lemma3}
There exists a class $\mathcal{KL}$ function $\beta_w$ such that for all solutions of the system satisfying (\ref{condition strictn1}) the following holds
\begin{eqnarray}
\Omega_w(t)&\leq&\beta_w\left(\Omega_w(0),t\right),\quad\mbox{for all $t\geq0$}\label{estimatew}\\
\Omega_w(t)&=&|X(t)|+\|w(t)\|_{\infty}+\|w_x(t)\|_{\infty}.\label{Omw}
\end{eqnarray}
\end{lemma}

In Lemmas \ref{lemma4}--\ref{lemma6} below, the equivalency of the $C^1$ norm between the original state variables $\left(X,u(x)\right)$, $x\in[0,1]$, and the state variables $\left(X,w(x)\right)$, $x\in[0,1]$, is established. The proofs of each of Lemmas \ref{lemma4} and \ref{lemma5}, utilize different arguments and employ different assumptions. For this reason, the proof of the norm equivalency, between the original and the new state variables, is decomposed into three different lemmas.


\begin{lemma}
\label{lemma4}
There exists a class $\mathcal{K}_{\infty}$ function $\rho_1$ such that for all solutions of the system satisfying (\ref{condition strictn1}) the following holds 
\begin{eqnarray}
\|p(t)\|_{\infty}+\|{p}_x(t)\|_{\infty}&\leq&\rho_1\left(|X(t)|+\|u(t)\|_{\infty}\right),\nonumber\\
&& \mbox{for all $t\geq0$}.\label{estimate f p}
\end{eqnarray}
\end{lemma}

\begin{lemma}
\label{lemma5}
There exists a class $\mathcal{K}_{\infty}$ function $\rho_2$ such that for all solutions of the system satisfying (\ref{condition strictn1}) the following holds 
\begin{eqnarray}
\|p(t)\|_{\infty}+\|{p}_x(t)\|_{\infty}&\leq&\rho_2\left(|X(t)|+\|w(t)\|_{\infty}\right),\nonumber\\
&& \mbox{for all $t\geq0$}.\label{estimate f p w}
\end{eqnarray}
\end{lemma}


\begin{lemma}
\label{lemma6}
There exist class $\mathcal{K}_{\infty}$ functions $\rho_3$ and $\rho_4$ such that for all solutions of the system satisfying (\ref{condition strictn1}) the following hold
\begin{eqnarray}
\Omega_w(t)&\leq& \rho_3\left(\Omega(t)\right),\quad \mbox{for all $t\geq0$}\label{est0w}\\
\Omega(t)&\leq& \rho_4\left(\Omega_w(t)\right),\quad \mbox{for all $t\geq0$},\label{est1w}
\end{eqnarray}
where $\Omega_w$ is defined in (\ref{Omw}) and $\Omega$ is defined in (\ref{omeg}). 
\end{lemma}

An estimate of the region of attraction of the predictor-feedback control law (\ref{controller})--(\ref{gamma}) within the feasibility region, defined by condition (\ref{condition strictn1}), is derived in the next two lemmas.

\begin{lemma}
\label{lemma7}
There exists a positive constant $\delta_1$ such that all of the solutions that satisfy
\begin{eqnarray}
|X(t)|+\|u(t)\|_{\infty}+\|u_x(t)\|_{\infty}< \delta_1,\quad \mbox{for all $t\geq0$},\label{gg}
\end{eqnarray}
they also satisfy (\ref{condition strictn1}).
\end{lemma}

\begin{lemma}
\label{lemma8}
There exists a positive constant $\delta$ such that for all initial conditions of the closed-loop system (\ref{pdesys1})--(\ref{pdesys2}), (\ref{controller})--(\ref{gamma}) that satisfy (\ref{region}), the solutions of the system satisfy (\ref{gg}), and hence, satisfy (\ref{condition strictn1}).
\end{lemma}

\paragraph*{Proof of Theorem \ref{thm1}}
Estimate (\ref{estimate}) in Theorem \ref{thm1} is proved combining Lemmas \ref{lemma3} and \ref{lemma6} with 
\begin{eqnarray}
\beta_u(s,t)=\rho_4\left(\beta_w\left(\rho_3\left(s\right),t\right)\right).
\end{eqnarray}

We show next the well-posedness of the system. We start by proving the well-posedness of the predictor 
\begin{eqnarray}
P(t)=p(1,t), 
\end{eqnarray}
where $p$ is defined in (\ref{p1}). Differentiating definition (\ref{p1}) with respect to $t$ and using integration by parts in the integral,  taking into account that $p$ satisfies $p_t(x,t)=v\left(u(x,t)\right)p_x(x,t)$ (see relation (\ref{first}) in Appendix~A) and employing relations (\ref{pdesys2}), (\ref{controller}) we obtain that
\begin{eqnarray}
\dot{P}(t)&=&v\left(\kappa\left(P(t)\right)\right)f\left(P(t),\kappa\left(P(t)\right)\right)\nonumber\\
&&\times\Gamma\left(u(1,t),u_x(1,t),1\right).\label{pop1}
\end{eqnarray}
From the definition of $\Gamma$ in (\ref{p1}) and (\ref{gamma}), using (\ref{pdesys3}) we obtain from (\ref{pdesys2}), (\ref{controller}) that
\begin{eqnarray}
\Gamma\left(u(1,t),u_x(1,t),1\right)&=&\frac{1}{v\left(\kappa\left(P(t)\right)\right)}\nonumber\\
&&-\frac{v'\left(\kappa\left(P(t)\right)\right)u_t(1,t)}{v\left(\kappa\left(P(t)\right)\right)^3},
\end{eqnarray}
and thus, from (\ref{pop1}) we arrive at
\begin{eqnarray}
\frac{1}{v\left(\kappa\left(P(t)\right)\right)}&=&\left(1+\frac{v'\left(\kappa\left(P(t)\right)\right)}{v\left(\kappa\left(P(t)\right)\right)^2}\right.\nonumber\\
&&\times\left.\vphantom{\frac{v'\left(\kappa\left(P(t)\right)\right)}{v\left(\kappa\left(P(t)\right)\right)^2}}\frac{\partial \kappa\left(P(t)\right)}{\partial P}f\left(P(t),\kappa\left(P(t)\right)\right)\right)\nonumber\\
&&\times\Gamma\left(u(1,t),u_x(1,t),1\right).\label{forg}
\end{eqnarray}
Since condition (\ref{condition strictn1}) guarantees the positivity and boundness of $\Gamma$ (see also relations (\ref{pos g}), (\ref{mu}) in Appendix~A), it follows from (\ref{forg}) and Assumption \ref{ass D} that the term $1+\frac{v'\left(\kappa\left(P(t)\right)\right)\frac{\partial \kappa\left(P(t)\right)}{\partial P}f\left(P(t),\kappa\left(P(t)\right)\right)}{v\left(\kappa\left(P(t)\right)\right)^2}$ is positive. Hence, solving (\ref{forg}) with respect to $\Gamma$ and substituting the resulting expression into (\ref{pop1}) we arrive at
\begin{eqnarray}
\dot{P}(t)&=&\frac{1}{1+\frac{v'\left(\kappa\left(P(t)\right)\right)\frac{\partial \kappa\left(P(t)\right)}{\partial P}f\left(P(t),\kappa\left(P(t)\right)\right)}{v\left(\kappa\left(P(t)\right)\right)^2}}\nonumber\\
&&\times f\left(P(t),\kappa\left(P(t)\right)\right).\label{popex}
\end{eqnarray}
Under the regularity assumptions on $v$ and $\kappa$, which follow from Assumptions \ref{ass D} and \ref{assi}, respectively, one can conclude that the right-hand side of (\ref{popex}) is Lipschitz with respect to $P$. Therefore, the compatibility conditions (\ref{comb1}), (\ref{comb2}) guarantee that there exists a unique solution $P(t)\in C^1[0,\infty)$. Moreover, employing estimates (\ref{estimate}), (\ref{estimate f p}), one can conclude that $\delta$ in the statement of Theorem \ref{thm1} can be chosen sufficiently small such that all the conditions of Theorem 1.1 in \cite{TLI} (Chapter 5) are satisfied, and hence, the existence and uniqueness of $u(x,t)\in C^1\left([0,1]\times [0,\infty)\right)$, which satisfies (\ref{pdesys3}), (\ref{pdesys2}), (\ref{controller}), follows (see also the discussion in, e.g., \cite{coron}, \cite{prier}). The fact that $u(x,t)\in C^1\left([0,1]\times [0,\infty)\right)$ and the regularity properties of $f$ imply from (\ref{pdesys1}) the existence and uniqueness of $X(t)\in C^1[0,\infty)$.

\section{Conclusions and Future Work}
\label{sec con}
We presented a predictor-feedback control design methodology for nonlinear systems with actuator dynamics governed by quasilinear, first-order hyperbolic PDEs. We proved that the closed-loop system, under the developed feedback law, is locally asymptotically stable, utilizing Lyapunov-like arguments and ISS estimates. We also emphasized the relation of the considered PDE-ODE cascade to a system with input delay that depends on past input values. 

A topic of future research may be the problem of boundary stabilization of general, quasilinear systems of first-order hyperbolic PDEs coupled with nonlinear ODE systems, as it is done in \cite{di meglio1} for the case in which both the PDE and ODE parts of the system are linear.


%
\setcounter{equation}{0}
\renewcommand{\theequation}{A.\arabic{equation}}
\section*{Appendix A}
\subsection*{Proof of Lemma \ref{lemma2}}
We first show that 
\begin{eqnarray}
p_t(x,t)=v\left(u(x,t)\right)p_x(x,t).\label{first}
\end{eqnarray}
Differentiating (\ref{p1}) with respect to $t$ and using (\ref{pdesys1}), (\ref{pdesys3}) as well as the fact that $p(0,t)=X(t)$, which immediately follows from (\ref{p1}) with $x=0$, we get that
\begin{eqnarray}
p_t(x,t)&=&f\left(p(0,t),u(0,t)\right)+\int_0^x\frac{\partial f\left(p(y,t),u(y,t)\right)}{\partial p}\nonumber\\
&&\times p_t(y,t)\Gamma\left(u(y,t),u_y(y,t),y\right)dy\nonumber\\
&&+\int_0^x\frac{\partial f\left(p(y,t),u(y,t)\right)}{\partial u}v\left(u(y,t)\right)u_y(y,t)\nonumber\\
&&\times\Gamma\left(u(y,t),u_y(y,t),y\right)dy\nonumber\\
&&+\int_0^x f\left(p(y,t),u(y,t)\right)\Gamma_u\left(u(y,t),u_y(y,t),y\right)\nonumber\\
&&\times v\left(u(y,t)\right)u_y(y,t)dy+\int_0^xf\left(p(y,t),u(y,t)\right)\nonumber\\
&&\times\Gamma_{u_y}\left(u(y,t),u_y(y,t),y\right)\left(v'\left(u(y,t)\right)u_y(y,t)^2\right.\nonumber\\
&&\left.\vphantom{v'\left(u(y,t)\right)}+v\left(u(y,t)\right)u_{yy}(y,t)\right)dy.\label{a1}
\end{eqnarray}
Differentiating (\ref{p1}) with respect to $x$ we get that
\begin{eqnarray}
v\left(u(x,t)\right)p_x(x,t)&=&\int_0^x\frac{\partial \Lambda(y,t)}{\partial y}dy+v\left(u(0,t)\right)\nonumber\\
&&\times f\left(p(0,t),u(0,t)\right)\nonumber\\
&&\times\Gamma\left(u(0,t),u_y(0,t),0\right)\\
\Lambda(y,t)&=&v\left(u(y,t)\right)f\left(p(y,t),u(y,t)\right)\nonumber\\
&&\times\Gamma\left(u(y,t),u_y(y,t),y\right),
\end{eqnarray}
and hence,
\begin{eqnarray}
v\left(u(x,t)\right)p_x(x,t)&=&v\left(u(0,t)\right)f\left(p(0,t),u(0,t)\right)\nonumber\\
&&\times\Gamma\left(u(0,t),u_y(0,t),0\right)\nonumber\\
&&+\int_0^x\frac{\partial f\left(p(y,t),u(y,t)\right)}{\partial p}\nonumber\\
&&\times p_y(y,t)v\left(u(y,t)\right)\nonumber\\
&&\times\Gamma\left(u(y,t),u_y(y,t),y\right)dy\nonumber\\
&&+\int_0^x\frac{\partial f\left(p(y,t),u(y,t)\right)}{\partial u}\nonumber\\
&&\times u_y(y,t)v\left(u(y,t)\right)\nonumber\\
&&\times\Gamma\left(u(y,t),u_y(y,t),y\right)dy\nonumber\\
&&+\int_0^x f\left(p(y,t),u(y,t)\right)\nonumber\\
&&\times\Gamma_u\left(u(y,t),u_y(y,t),y\right)\nonumber\\
&&\times v\left(u(y,t)\right) u_y(y,t)dy\nonumber\\
&&+\int_0^x f\left(p(y,t),u(y,t)\right)\nonumber\\
&&\times\Gamma_{u_y}\left(u(y,t),u_y(y,t),y\right)\nonumber\\
&&\times v\left(u(y,t)\right) u_{yy}(y,t)dy\nonumber\\
&&+\int_0^x f\left(p(y,t),u(y,t)\right)\nonumber\\
&&\times\Gamma_y\left(u(y,t),u_y(y,t),y\right)\nonumber\\
&&\times v\left(u(y,t)\right)dy\nonumber\\
&&+\int_0^x f\left(p(y,t),u(y,t)\right)\nonumber\\
&&\times\Gamma\left(u(y,t),u_y(y,t),y\right)\nonumber\\
&&\times v'\left(u(y,t)\right) u_y(y,t)dy.\label{a4}
\end{eqnarray}
Comparing (\ref{a1}) with (\ref{a4}) and using the fact that $\Gamma\left(u(0,t),u_y(0,t),0\right)=\frac{1}{v\left(u(0,t)\right)}$, which follows from (\ref{gamma}) for $y=0$, we arrive at
\begin{eqnarray}
P(x,t)&=&\int_0^x\frac{\partial f\left(p(y,t),u(y,t)\right)}{\partial p}P(y,t)\nonumber\\
&&\times\Gamma\left(u(y,t),u_y(y,t),y\right)dy\nonumber\\
&&+\int_0^xf\left(p(y,t),u(y,t)\right)\nonumber\\
&&\times\Gamma_{u_y}\left(u(y,t),u_y(y,t),y\right)\nonumber\\
&&\times v'\left(u(y,t)\right)u_y(y,t)^2dy\nonumber\\
&&-\int_0^x f\left(p(y,t),u(y,t)\right)\nonumber\\
&&\times\Gamma_y\left(u(y,t),u_y(y,t),y\right)\nonumber\\
&&\times v\left(u(y,t)\right)dy\nonumber\\
&&-\int_0^x f\left(p(y,t),u(y,t)\right)\nonumber\\
&&\times\Gamma\left(u(y,t),u_y(y,t),y\right)\nonumber\\
&&\times v'\left(u(y,t)\right)u_y(y,t)dy,\label{po}
\end{eqnarray}
where we defined
\begin{eqnarray}
P(x,t)=p_t(x,t)-v\left(u(x,t)\right)p_x(x,t).
\end{eqnarray}
Using the definition of $\Gamma$ in (\ref{gamma}) we get that
\begin{eqnarray}
\Gamma_{u_y}\left(u(y,t),u_y(y,t),y\right)&=&-\frac{yv'\left(u(y,t)\right)}{v\left(u(y,t)\right)^2}\label{gas}\\
\Gamma_{y}\left(u(y,t),u_y(y,t),y\right)&=&-\frac{v'\left(u(y,t)\right)u_y(y,t)}{v\left(u(y,t)\right)^2}\label{gas1}.
\end{eqnarray}
Combining (\ref{gas}), (\ref{gas1}) and using (\ref{gamma}) we arrive at
\begin{eqnarray}
\Gamma_{u_y}\left(u,u_y,y\right)v'\left(u\right)u_y^2&=&\Gamma_y\left(u,u_y,y\right)v\left(u\right)\nonumber\\
&&-\frac{yv'\left(u\right)^2u_y^2}{v\left(u\right)^2}+\frac{v'\left(u\right)u_y}{v\left(u\right)}\label{equal1}\\
\Gamma\left(u,u_y,y\right)v'\left(u\right)u_y&=&\frac{v'\left(u\right)u_y}{v\left(u\right)}-\frac{yv'\left(u\right)^2u_y^2}{v\left(u\right)^2}.\label{equal2}
\end{eqnarray}
Since the right-hand sides of (\ref{equal1}), (\ref{equal2}) are equal, from (\ref{po}) it follows that
\begin{eqnarray}
P(x,t)&=&\int_0^x\frac{\partial f\left(p(y,t),u(y,t)\right)}{\partial p}P(y,t)\nonumber\\
&&\times\Gamma\left(u(y,t),u_y(y,t),y\right)dy.\label{po1}
\end{eqnarray}
Therefore, for each $t\geq0$, the function $P$ satisfies for all $x\in[0,1]$
\begin{eqnarray}
P_x(x,t)&=&\frac{\partial f\left(p(x,t),u(x,t)\right)}{\partial p}\nonumber\\
&&\times\Gamma\left(u(x,t),u_x(x,t),x\right)P(x,t)\\
P(0,t)&=&0.
\end{eqnarray}
Hence, 
\begin{eqnarray}
P&\equiv&0,
\end{eqnarray}
which proves that indeed (\ref{first}) holds. Therefore, differentiating (\ref{w1}) with respect to $t$ we get that
\begin{eqnarray}
w_t(x,t)&=&u_t(x,t)+\frac{\partial \kappa\left(p(x,t)\right)}{\partial p}v\left(u(x,t)\right)\nonumber\\
&&\times p_x(x,t).\label{we}
\end{eqnarray}
Differentiating (\ref{w1}) with respect to $x$ we get that
\begin{eqnarray}
v\left(u(x,t)\right)w_x(x,t)&=&v\left(u(x,t)\right)u_x(x,t)\nonumber\\
&&+\frac{\partial \kappa\left(p(x,t)\right)}{\partial p}v\left(u(x,t)\right)\nonumber\\
&&\times p_x(x,t).\label{we1}
\end{eqnarray}
Combining (\ref{we}) with (\ref{we1}) and using (\ref{pdesys3}) we arrive at (\ref{pdent}). Furthermore, since from (\ref{p1}) it holds that $p(0,t)=X(t)$, relation (\ref{pdesysn1t}) follows from (\ref{pdesys1}) and (\ref{w1}) for $x=0$. Finally, relation (\ref{pden1t}) follows from (\ref{w1}) for $x=1$ and (\ref{controller}), (\ref{pdesys2}).

\subsection*{Proof of Lemma \ref{lemma3}}
Consider the following Lyapunov functional
\begin{eqnarray}
L_{c,m}(t)&=&\int_0^1e^{2(c+\lambda)xm}w(x,t)^{2m}dx\nonumber\\
&&+\int_0^1e^{2(c+\lambda)xm}w_x(x,t)^{2m}dx,\label{L1}
\end{eqnarray}
for any $c>0$ and any positive integer $m$, where (under Assumption \ref{ass D})
\begin{eqnarray}
w_{xt}(x,t)&=&v'\left(u(x,t)\right)u_x(x,t)w_x(x,t)\nonumber\\
&&+v\left(u(x,t)\right)w_{xx}(x,t)\label{pdentnew}\\
w_x(1,t)&=&0.\label{pden1tnew}
\end{eqnarray}
Under Assumption \ref{ass D} (positivity of $v$), taking the time derivative of (\ref{L1}) along the solutions of (\ref{pdent})--(\ref{pdesysn1t}), (\ref{pdentnew}), (\ref{pden1tnew}) we get using integration by parts that
\begin{eqnarray}
\dot{L}_{c,m}(t)&\leq&-\int_0^1e^{2(c+\lambda)xm}w(x,t)^{2m}\nonumber\\
&&\times\left(2m(c+\lambda)v\left(u(x,t)\right)\right.\nonumber\\
&&\left.+v'\left(u(x,t)\right)u_x(x,t)\right)dx\nonumber\\
&&-\int_0^1e^{2(c+\lambda)xm}w_x(x,t)^{2m}\nonumber\\
&&\times\left(2m(c+\lambda)v\left(u(x,t)\right)+v'\left(u(x,t)\right)u_x(x,t)\right.\nonumber\\
&&\left.-2mv'\left(u(x,t)\right)u_x(x,t)\right)dx.\label{inter1}
\end{eqnarray}
Using (\ref{condition strictn1}) and the fact that $m\geq1$ we obtain for all $x\in[0,1]$ and $t\geq0$
\begin{eqnarray}
2m\left(-\lambda+M\right)v\left(u(x,t)\right)&\geq&-2m\lambda v\left(u(x,t)\right)\nonumber\\
&&-v'\left(u(x,t)\right)u_x(x,t)\\
2m\left(-\lambda+1\right)v\left(u(x,t)\right)&\geq&-2m\lambda v\left(u(x,t)\right)-(1-2m)\nonumber\\
&&\times v'\left(u(x,t)\right)u_x(x,t).
\end{eqnarray}
Therefore, choosing any $\lambda$ such that $\lambda\geq 1+M$, it follows from (\ref{inter1}) that
\begin{eqnarray}
\dot{L}_{c,m}(t)&\leq&-2mc\underline{v}\int_0^1e^{2(c+\lambda)xm}w(x,t)^{2m}\nonumber\\
&&-2mc\underline{v}\int_0^1e^{2(c+\lambda)xm}w_x(z,t)^{2m}dx,
\end{eqnarray}
where we also used (\ref{assv}). Thus\footnote{Note that although the estimate (\ref{distribution}) for $L_{c,m}$ is derived for $u$ that is of class $C^2$ (and thus, so is $w$ satisfying (\ref{pdentnew}), (\ref{pden1tnew})), the estimate (\ref{distribution}) remains valid (in the distribution sense) when $u$ is only of class $C^1$, see, e.g., \cite{coron}. },
\begin{eqnarray}
\dot{L}_{c,m}(t)&\leq&-2mc\underline{v} L_{c,m}(t),\label{distribution}
\end{eqnarray}
which implies that
\begin{eqnarray}
{L}^{\frac{1}{2m}}_{c,m}(t)&\leq e^{-c\underline{v} (t-s)} L^{\frac{1}{2m}}_{c,m}(s),\quad\mbox{for all $t\geq s\geq0$}.
\end{eqnarray}
Moreover, from (\ref{L1}) it follows that
\begin{eqnarray}
\Xi_{c,m}(t)\leq  2e^{-c\underline{v} (t-s)}\Xi_{c,m}(s),\quad\mbox{for all $t\geq s\geq0$},\label{takel}
\end{eqnarray}
where
\setlength{\arraycolsep}{5pt}\begin{eqnarray}
\Xi_{c,m}(t)&=&\left(\int_0^1e^{2(c+\lambda)xm}w(x,t)^{2m}dx\right)^{\frac{1}{2m}}\nonumber\\
&&+\left(\int_0^1e^{2(c+\lambda)xm}w_x(x,t)^{2m}dx\right)^{\frac{1}{2m}}.\label{definition XI}
\end{eqnarray}\setlength{\arraycolsep}{5pt}Taking the limit of (\ref{takel}) as $m$ goes to infinity, with the definition of the maximum norm, i.e., with relation $\|\theta(t)\|_{\infty}=\lim_{m\to\infty}\left(\int_0^1|\theta(x,t)|^{2m}dx\right)^{\frac{1}{2m}}$, we obtain from (\ref{definition XI})
\begin{eqnarray}
\Xi_c(t)\leq2e^{-c\underline{v} (t-s)}\Xi_c(s),\quad\mbox{for all $t\geq s\geq0$}.
\end{eqnarray}
where
\begin{eqnarray}
\Xi_c(t)&=&\max_{0\leq x\leq 1}\left|e^{x(c+\lambda)}w(x,t)\right|\nonumber\\
&&+\max_{0\leq x\leq 1}\left|e^{x(c+\lambda)}w_x(x,t)\right|.
\end{eqnarray}
It follows, for all $t\geq s\geq0$, that
\setlength{\arraycolsep}{5pt}\begin{eqnarray}
\|w(t)\|_{\infty}+\|w_x(t)\|_{\infty}&\leq& 2e^{-c\underline{v} (t-s)}e^{(c+\lambda)}\nonumber\\
&&\!\!\times\left(\|w(s)\|_{\infty}+\|w_x(s)\|_{\infty} \right)\!.\label{westimaten}
\end{eqnarray}\setlength{\arraycolsep}{5pt}
Under Assumption \ref{assi} (see, e.g., \cite{sontag}) we obtain from (\ref{pdesysn1t}) that
\begin{eqnarray}
|X(t)|\leq\beta_1\left(|X(s)|,t-s\right)+\gamma_1\left(\sup_{s\leq \tau\leq t}\left|w(0,\tau)\right|\right)\!,\label{kh1}
\end{eqnarray}
for all $t\geq s\geq0$, some class $\mathcal{KL}$ function $\beta_1$, and some class $\mathcal{K}$ function $\gamma_1$. Mimicking the arguments in the proof of Lemma 4.7 from \cite{khalil}, we set $s=\frac{t}{2}$ in (\ref{kh1}) to get that
\begin{eqnarray}
|X(t)|&\leq&\beta_1\left(\left|X\left(\frac{t}{2}\right)\right|,\frac{t}{2}\right)\nonumber\\
&&+\gamma_1\left(\sup_{\frac{t}{2}\leq \tau\leq t}\|w(\tau)\|_{\infty}\right),\label{kh2}
\end{eqnarray}
and thus, using (\ref{kh1}) for $s=0$ and $t\to \frac{t}{2}$ we arrive at
\setlength{\arraycolsep}{5pt}\begin{eqnarray}
|X(t)|&\leq&\beta_1\left(\vphantom{\gamma_1\left(\sup_{0\leq \tau\leq \frac{t}{2}}\|w(\tau)\|_{\infty}\right)}\beta_1\left(|X(0)|,\frac{t}{2}\right)\right.\nonumber\\
&&\left.+\gamma_1\left(\sup_{0\leq \tau\leq \frac{t}{2}}\|w(\tau)\|_{\infty}\right),\frac{t}{2}\right)\nonumber\\
&&+\gamma_1\left(\sup_{\frac{t}{2}\leq \tau\leq t}\|w(\tau)\|_{\infty}\right),\!\!\!\quad\mbox{for all $t\geq0$}.\label{kh3}
\end{eqnarray}\setlength{\arraycolsep}{5pt}Moreover, using (\ref{westimaten}) we get that
\setlength{\arraycolsep}{5pt}\begin{eqnarray}
\sup_{0\leq \tau\leq \frac{t}{2}}\|w(\tau)\|_{\infty}&\leq&2e^{(c+\lambda)} \left(\|w(0)\|_{\infty}+\|w_x(0)\|_{\infty} \right)\label{kh4}\\
\sup_{\frac{t}{2}\leq \tau\leq t}\|w(\tau)\|_{\infty}&\leq&2e^{-c\underline{v} \frac{t}{2}}e^{(c+\lambda)} \left(\|w(0)\|_{\infty}\right.\nonumber\\
&&\left.+\|w_x(0)\|_{\infty} \right).\label{kh5}
\end{eqnarray}\setlength{\arraycolsep}{5pt}Therefore, combining (\ref{kh3}) with (\ref{kh4}), (\ref{kh5}) and using (\ref{westimaten}) we get (\ref{estimatew}) with
\setlength{\arraycolsep}{0pt}\begin{eqnarray}
\beta_w\left(s,t\right)&=&\beta_1\left(\beta_1\left(s,0\right)+\gamma_1\left(2e^{(c+\lambda)} s\right),\frac{t}{2}\right)\nonumber\\
&&+2e^{-c\underline{v} t}e^{(c+\lambda)}s+\gamma_1\left(2e^{-c\underline{v} \frac{t}{2}}e^{(c+\lambda)}s\right) .\label{kh6}
\end{eqnarray}\setlength{\arraycolsep}{5pt}

\subsection*{Proof of Lemma \ref{lemma4}}
Differentiating relation (\ref{p1}) with respect to $x$ we get that, for each $t$, $p$ satisfies the following ODE in $x$
\begin{eqnarray}
p_x(x,t)&=&f\left(p(x,t),u(x,t)\right)\Gamma\left(u(x,t),u_x(x,t),x\right)\label{odep}\\
p(0,t)&=&X(t).\label{odep1}
\end{eqnarray}
Under Assumption \ref{assf} there exists a smooth function $R:\mathbb{R}^n\to\mathbb{R}_+$ and class $\mathcal{K}_{\infty}$ functions $\alpha_1$, $\alpha_2$, and $\alpha_3$ such that (see, e.g., \cite{angeli}, \cite{krstic book}, \cite{krstic})
\begin{eqnarray}
\alpha_1\left(|X|\right)&\leq&R\left(X\right)\leq\alpha_2\left(|X|\right)\label{a9}\\
\frac{\partial R\left(X\right)}{\partial X}f\left(X,\omega\right)&\leq& R\left(X\right)+\alpha_3\left(|\omega|\right),\label{a8}
\end{eqnarray}
for all $\left(X,\omega\right)^{\rm T}\in\mathbb{R}^{n+1}$. From relation (\ref{condition strictn1}) it follows that
\begin{eqnarray}
v\left(u(x,t)\right)-xv'\left(u(x,t)\right)u_x(x,t)&>&0,\quad \mbox{for all $x\in[0,1]$}\nonumber\\
&&\!\!\qquad\mbox{ and $t\geq0$},
\end{eqnarray}
which can be seen considering separately the cases $v'\left(u(x,t)\right)u_x(x,t)\leq0$ and $v'\left(u(x,t)\right)u_x(x,t)>0$, and using (\ref{assv}). Hence, from the definition of $\Gamma$ in (\ref{gamma}) and (\ref{assv}) we conclude that 
\begin{eqnarray}
\Gamma\left(u(x,t),u_x(x,t),x\right)&>&0,\nonumber\\
&& \mbox{for all $x\in[0,1]$ and $t\geq0$}. \label{pos g}
\end{eqnarray}
Thus, from (\ref{a8}) it follows that
\begin{eqnarray}
\Pi(x,t)&\leq& \Gamma\left(u(x,t),u_x(x,t),x\right)\left(R\left(p(x,t)\right)\right.\nonumber\\
&&+\left.\alpha_3\left(|u(x,t)|\right)\right)\label{comp1}\\
\Pi(x,t)&=&\frac{\partial R\left(p(x,t)\right)}{\partial p}f\left(p(x,t),u(x,t)\right)\nonumber\\
&&\times \Gamma\left(u(x,t),u_x(x,t),x\right).
\end{eqnarray}
From (\ref{gamma}), using relations (\ref{assv}) and (\ref{condition strictn1}) it also follows that
\begin{eqnarray}
\Gamma\left(u(x,t),u_x(x,t),x\right)&\leq&\frac{2+M}{\underline{v}},\nonumber\\
&& \mbox{for all $x\in[0,1]$ and $t\geq0$},\label{mu}
\end{eqnarray}
and thus, from (\ref{comp1}) we get using (\ref{odep}) that
\begin{eqnarray}
\frac{\partial R\left(p(x,t)\right)}{\partial x}&\!\leq\!& \frac{2+M}{\underline{v}}\left(R\left(p(x,t)\right)+\alpha_3\left(|u(x,t)|\right)\right).\label{comp}
\end{eqnarray}
Employing the comparison principle and using (\ref{odep1}) we arrive at
\begin{eqnarray}
R\left(p(x,t)\right)&\leq& e^{\frac{2+M}{\underline{v}} x}R\left(X(t)\right)+\frac{2+M}{\underline{v}}\nonumber\\
&&\times \int_0^x e^{\frac{2+M}{\underline{v}} (x-y)}\alpha_3\left(|u(y,t)|\right)dy,
\end{eqnarray}
for all $x\in[0,1]$ and $t\geq0$. Hence, using (\ref{a9}) we get
\begin{eqnarray}
\|p(t)\|_{\infty}\leq\alpha_4\left(|X(t)|+\|u(t)\|_{\infty}\right),\label{pop}
\end{eqnarray}
where 
\begin{eqnarray}
\alpha_4(s)=\alpha^{-1}_1\left(e^{\frac{2+M}{\underline{v}}}\left(\alpha_2(s)+\alpha_3(s)\right)\right).
\end{eqnarray}
Since $f$ is continuously differentiable with $f(0,0)=0$ we conclude that there exists a class $\mathcal{K}_{\infty}$ function $\alpha_5$ such that
\begin{eqnarray}
\left|f(X,\omega)\right|\leq\alpha_5\left(|X|+|\omega|\right).\label{ff}
\end{eqnarray}
Thus, using (\ref{pop}) and (\ref{mu}), we get from (\ref{odep})
\setlength{\arraycolsep}{5pt}\begin{eqnarray}
\left|{p}_x(x,t)\right|&\leq&\alpha_6\left(|X(t)|+\|u(t)\|_{\infty}\right),\label{both}
\end{eqnarray}\setlength{\arraycolsep}{5pt}where
\begin{eqnarray}
\alpha_6(s)=\frac{2+M}{\underline{v}}\alpha_5\left(\alpha_4(s)+s\right).
\end{eqnarray}
The proof is completed by taking the maximum, with respect to $x\in[0,1]$, in both sides of (\ref{both}) and setting $\rho_1(s)=\alpha_4(s)+\alpha_6(s)$.

\subsection*{Proof of Lemma \ref{lemma5}}
We start defining the change of variables with respect to $x$
\begin{eqnarray}
z(x,t)=\frac{x}{v\left(u(x,t)\right)},
\end{eqnarray}
which is well-defined thanks to (\ref{assv}) and where $t$ acts as a parameter. Using the fact that 
\begin{eqnarray}
\Gamma\left(u(x,t),u_x(x,t),x\right)=\frac{\partial z(x,t)}{\partial x},\label{def z}
\end{eqnarray}
it follows from (\ref{pos g}) that the function $z$ is strictly increasing with respect to $x$, for each $t$. Thus, it admits an inverse defined for each $t$ as $x=\chi(z,t)$. Therefore, from relations (\ref{odep}), (\ref{odep1}), and definition (\ref{def z}) we obtain
\begin{eqnarray}
 \bar{p}_z\left(z,t\right)&=&f\left(\bar{p}(z,t),\bar{u}(z,t)\right),\quad z\in\left[0,\frac{1}{v\left(u(1,t)\right)}\right]\label{fgg}\\
\bar{p}(0,t)&=&X(t),\label{fgg1}
\end{eqnarray}
where 
\begin{eqnarray}
\bar{p}(z,t)&=&p\left(\chi(z,t),t\right)\label{p1de}\\
\bar{u}(z,t)&=&u\left(\chi(z,t),t\right). \label{p2de}
\end{eqnarray}
Moreover, setting $x=\chi(z,t)$ in relation (\ref{w1}) we get that
\begin{eqnarray}
\bar{u}(z,t)=\kappa\left(\bar{p}(z,t)\right)+\bar{w}(z,t),
\end{eqnarray}
where 
\begin{eqnarray}
\bar{w}(z,t)=w\left(\chi(z,t),t\right). \label{p3de}
\end{eqnarray}
Thus, we re-write (\ref{fgg}) as 
\begin{eqnarray}
 \bar{p}_z\left(z,t\right)&=&f\left(\bar{p}(z,t),\kappa\left(\bar{p}(z,t)\right)+\bar{w}(z,t)\right),\nonumber\\
 && z\in\left[0,\frac{1}{v\left(u(1,t)\right)}\right].\label{fggnb}
\end{eqnarray}
Under Assumption \ref{assi} there exist a smooth function $S: \mathbb{R}^n\to\mathbb{R}_+$ and class $\mathcal{K}_{\infty}$ functions $\hat{\alpha}_1$, $\hat{\alpha}_2$, $\hat{\alpha}_3$, and $\hat{\alpha}_4$  such that (see, e.g., \cite{sontag})
\begin{eqnarray}
\hat{\alpha}_1\left(|X|\right)&\leq&S\left(X\right)\leq\hat{\alpha}_2\left(|X|\right)\label{a9iss}\\
\frac{\partial S\left(X\right)}{\partial X}f\left(X,\kappa\left(X\right)+\omega\right)&\leq& -\hat{\alpha}_3\left(|X|\right)+\hat{\alpha}_4\left(|\omega|\right).\label{a8iss}
\end{eqnarray}
Thus, we get from (\ref{a8iss}) that
\begin{eqnarray}
\Theta(z,t)&\leq& -\hat{\alpha}_3\left(\left|\bar{p}(z,t)\right|\right)+\hat{\alpha}_4\left(\left|\bar{w}(z,t)\right|\right)\label{kh1w}\\
\Theta(z,t)&=&\frac{\partial S\left(\bar{p}(z,t)\right)}{\partial \bar{p}}\nonumber\\
&&\times f\left(\bar{p}(z,t),\kappa\left(\bar{p}(z,t)\right)+\bar{w}(z,t)\right),
\end{eqnarray}
and hence, using (\ref{fggnb}), (\ref{fgg1}) and integrating from $0$ to $z$ we obtain
\begin{eqnarray}
S\left(\bar{p}(z,t)\right)&\leq& S\left(X(t)\right)+\int_0^z\hat{\alpha}_4\left(\left|\bar{w}(y,t)\right|\right)dy,\nonumber\\
&& z\in\left[0,\frac{1}{v\left(u(1,t)\right)}\right].\label{kh1wee}
\end{eqnarray}
Using (\ref{assv}) and (\ref{a9iss}) we get from (\ref{kh1wee}) that
\begin{eqnarray}
\left|\bar{p}(z,t)\right|&\leq&\hat{\alpha}_1^{-1}\left(\vphantom{+\frac{1}{\underline{v}}\hat{\alpha}_4\left(\max_{0\leq z\leq \frac{1}{v\left(u(1,t)\right)}}|\bar{w}(z,t)|\right)} \hat{\alpha}_2\left(\left|X(t)\right|\right)\right. \nonumber\\
&&\left.+\frac{1}{\underline{v}}\hat{\alpha}_4\left(\max_{0\leq z\leq \frac{1}{v\left(u(1,t)\right)}}|\bar{w}(z,t)|\right)\right),\nonumber\\
&&z\in\left[0,\frac{1}{v\left(u(1,t)\right)}\right].\label{kh1weenew}
\end{eqnarray}
Taking a maximum of both sides in (\ref{kh1weenew}), with definitions (\ref{p1de}), (\ref{p2de}), and (\ref{p3de}) we arrive at
\begin{eqnarray}
\|{p}(t)\|_{\infty}\leq\hat{\alpha}_5\left(\left|X(t)\right|+\|w(t)\|_{\infty}\right),\label{kh1wsd}
\end{eqnarray}
with $\hat{\alpha}_5(s)=\hat{\alpha}_1^{-1}\left(\hat{\alpha}_2(s)+\frac{1}{\underline{v}}\hat{\alpha}_4(s)\right)$. Under Assumption \ref{assi} (continuity of $\kappa$ and the fact that $\kappa(0)=0$) there exists a class $\mathcal{K}_{\infty}$ function $\hat{\bar{\alpha}}_1$ such that
\begin{eqnarray}
\left|\kappa\left(X\right)\right|&\leq& \hat{\bar{\alpha}}_1\left(|X|\right)\label{kk1}.
\end{eqnarray}
Therefore, using (\ref{ff}) we get from (\ref{fggnb}) that
\begin{eqnarray}
\left| \bar{p}_z\left(z,t\right)\right|&\!\leq\!&\alpha_5\left(\left|\bar{p}(z,t)\right|+\hat{\bar{\alpha}}_1\left(\left|\bar{p}(z,t)\right|\right)+\left|\bar{w}(z,t)\right|\right)\!.\label{done}
\end{eqnarray}
Using (\ref{kh1weenew}) we arrive at
\begin{eqnarray}
\left| \bar{p}_z\left(z,t\right)\right|\leq\hat{\alpha}_6\left(\left|X(t)\right|+\max_{0\leq z\leq \frac{1}{v\left(u(1,t)\right)}}|\bar{w}(z,t)|\right),\label{kh1wsdnew}
\end{eqnarray}
where $\hat{\alpha}_6(s)=\alpha_5\left(\hat{\alpha}_5(s)+\hat{\bar{\alpha}}_1\left(\hat{\alpha}_5(s)\right)+s\right)$. Using definition (\ref{p1de}), from relations (\ref{mu}), (\ref{def z}) we obtain that
\begin{eqnarray}
\left| {p}_x\left(x,t\right)\right|\leq\frac{2+M}{\underline{v}}\left| \bar{p}_z\left(z,t\right)\right|,
\end{eqnarray}
and hence, from (\ref{kh1wsdnew}) we arrive at
\begin{eqnarray}
\| {p}_x(t)\|_{\infty}&\leq&\frac{2+M}{\underline{v}}\hat{\alpha}_6\left(\vphantom{+\max_{0\leq z\leq \frac{1}{v\left(u(1,t)\right)}}|\bar{w}(z,t)|}\left|X(t)\right|\right.\nonumber\\
&&\left.+\max_{0\leq z\leq \frac{1}{v\left(u(1,t)\right)}}|\bar{w}(z,t)|\right).
\end{eqnarray}
With definition (\ref{p3de}) we get that
\begin{eqnarray}
\| {p}_x(t)\|_{\infty}\leq\frac{2+M}{\underline{v}}\hat{\alpha}_6\left(\left|X(t)\right|+\|{w}(t)\|_{\infty}\right),
\end{eqnarray}
and hence, the lemma is proved with $\rho_2(s)=\hat{\alpha}_5(s)+\frac{2+M}{\underline{v}}\hat{\alpha}_6\left(s\right)$.

\subsection*{Proof of Lemma \ref{lemma6}}
Under Assumption \ref{assi} (continuous differentiability of $\kappa$) there exists a class $\mathcal{K}_{\infty}$ function $\hat{\bar{\alpha}}_2$ such that
\begin{eqnarray}
\left|\nabla \kappa\left(X\right)\right|&\leq&\left|\nabla \kappa\left(0\right)\right|+ \hat{\bar{\alpha}}_2\left(|X|\right),\label{kk2}
\end{eqnarray}
for all $X\in\mathbb{R}^n$. Therefore, using (\ref{estimate f p}) and (\ref{kk1}), we get from (\ref{w1}) that
\setlength{\arraycolsep}{5pt}\begin{eqnarray}
\left|w(x,t)\right|+\left|w_x(x,t)\right|&\leq& |u(x,t)|+|u_x(x,t)|\nonumber\\
&&+\hat{\bar{\alpha}}_1\left(\rho_1\left(|X(t)|+\|u(t)\|_{\infty}\right)\right)\nonumber\\
&&+\left(\left|\nabla \kappa\left(0\right)\right|+\hat{\bar{\alpha}}_2\left(\rho_1\left(|X(t)|\right.\right.\right.\nonumber\\
&&\left.\left.\left.+\|u(t)\|_{\infty}\right)\right)\right)\nonumber\\
&&\times\rho_1\left(|X(t)|+\|u(t)\|_{\infty}\right),\label{fun}
\end{eqnarray}\setlength{\arraycolsep}{5pt}and thus, taking a maximum in both sides of (\ref{fun}), estimate (\ref{est0w}) follows with
\begin{eqnarray}
\rho_3\left(s\right)&=&s+\hat{\bar{\alpha}}_1\left(\rho_1\left(s\right)\right)+\left(\left|\nabla \kappa\left(0\right)\right|+\hat{\bar{\alpha}}_2\left(\rho_1\left(s\right)\right)\right)\nonumber\\
&&\times\rho_1\left(s\right).
\end{eqnarray}
Similarly, using (\ref{kk1}), (\ref{kk2}), and (\ref{estimate f p w}) we get estimate (\ref{est1w}) with
\begin{eqnarray}
\rho_4\left(s\right)&=&s+\hat{\bar{\alpha}}_1\left(\rho_2\left(s\right)\right)+\left(\left|\nabla \kappa\left(0\right)\right|+\hat{\bar{\alpha}}_2\left(\rho_2\left(s\right)\right)\right)\nonumber\\
&&\times\rho_2\left(s\right).
\end{eqnarray}

\subsection*{Proof of Lemma \ref{lemma7}}
Under Assumption \ref{ass D} (continuous differentiability of $v$) we conclude that there exists a class $\mathcal{K}_{\infty}$ function $\hat{\rho}$ such that
\begin{eqnarray}
\left|v'\left(u(x,t)\right)\right|&\leq&\left|v'(0)\right|+\hat{\rho}\left(\left|u(x,t)\right|\right),\label{sh1}
\end{eqnarray}
and hence, for all $x\in[0,1]$ and $t\geq0$ it holds that
\begin{eqnarray}
\left|v'\left(u(x,t)\right)\right|&\leq&\left|v'(0)\right|+\hat{\rho}\left(\|u(t)\|_{\infty}\right).\label{sh1new}
\end{eqnarray}
Thus, it holds that
\begin{eqnarray}
\left|v'\left(u(x,t)\right)u_x(x,t)\right|&\leq&\left(\left|v'(0)\right|+\hat{\rho}\left(\|u(t)\|_{\infty}\right)\right)\nonumber\\
&&\times\|u_x(t)\|_{\infty},\label{sh1newgh}
\end{eqnarray}
for all $x\in[0,1]$ and $t\geq0$. From relation (\ref{assv}), one can conclude that whenever
\begin{eqnarray}
\left|v'\left(u(x,t)\right)u_x(x,t)\right|&\leq&\epsilon,
\end{eqnarray}
where $\epsilon$ is any constant such that $0<\epsilon<\underline{v}$, relation (\ref{condition strictn1}) holds with any $M$ such that $M>\frac{\epsilon}{\underline{v}}$. Consequently, choosing any constant $\delta_1$ such that
\begin{eqnarray}
\delta_1\leq\psi^{-1}\left(\epsilon\right),
\end{eqnarray}
where 
\begin{eqnarray}
\psi\left(s\right)=\left(\left|v'(0)\right|+\hat{\rho}\left(s\right)\right)s,
\end{eqnarray}
completes the proof.

\subsection*{Proof of Lemma \ref{lemma8}}
Combining estimate (\ref{est1w}) with (\ref{estimatew}) we obtain
\begin{eqnarray}
\Omega(t)&\leq& \rho_4\left(\beta_w\left(\Omega_w(0),t\right)\right),
\end{eqnarray}
and hence, with (\ref{est0w}) and the properties of class $\mathcal{KL}$ functions we arrive at
\begin{eqnarray}
\Omega(t)&\leq& \rho_4\left(\beta_w\left(\rho_3\left(\Omega(0)\right),0\right)\right).
\end{eqnarray}
Therefore, for all initial conditions that satisfy the bound (\ref{region}) with any $\delta$ such that
\begin{eqnarray}
\delta\leq\phi^{-1}\left(\delta_1\right),
\end{eqnarray}
where
\begin{eqnarray}
\phi\left(s\right)= \rho_4\left(\beta_w\left(\rho_3\left(s,0\right)\right)\right),
\end{eqnarray}
the solutions satisfy (\ref{gg}). Furthermore, from Lemma \ref{lemma7}, it follows that for all of those initial conditions, the solutions verify (\ref{condition strictn1}).

\section*{Acknowledgments}
\textcolor{black}{Nikolaos Bekiaris-Liberis was supported by the funding from the European Commission's Horizon 2020 research and innovation programme under the Marie Sklodowska-Curie grant agreement No. 747898, project PADECOT.}

 {\small

\begin{wrapfigure}{L}{0.18\textwidth}
\vspace{-3mm}
\includegraphics[width=0.21\textwidth]{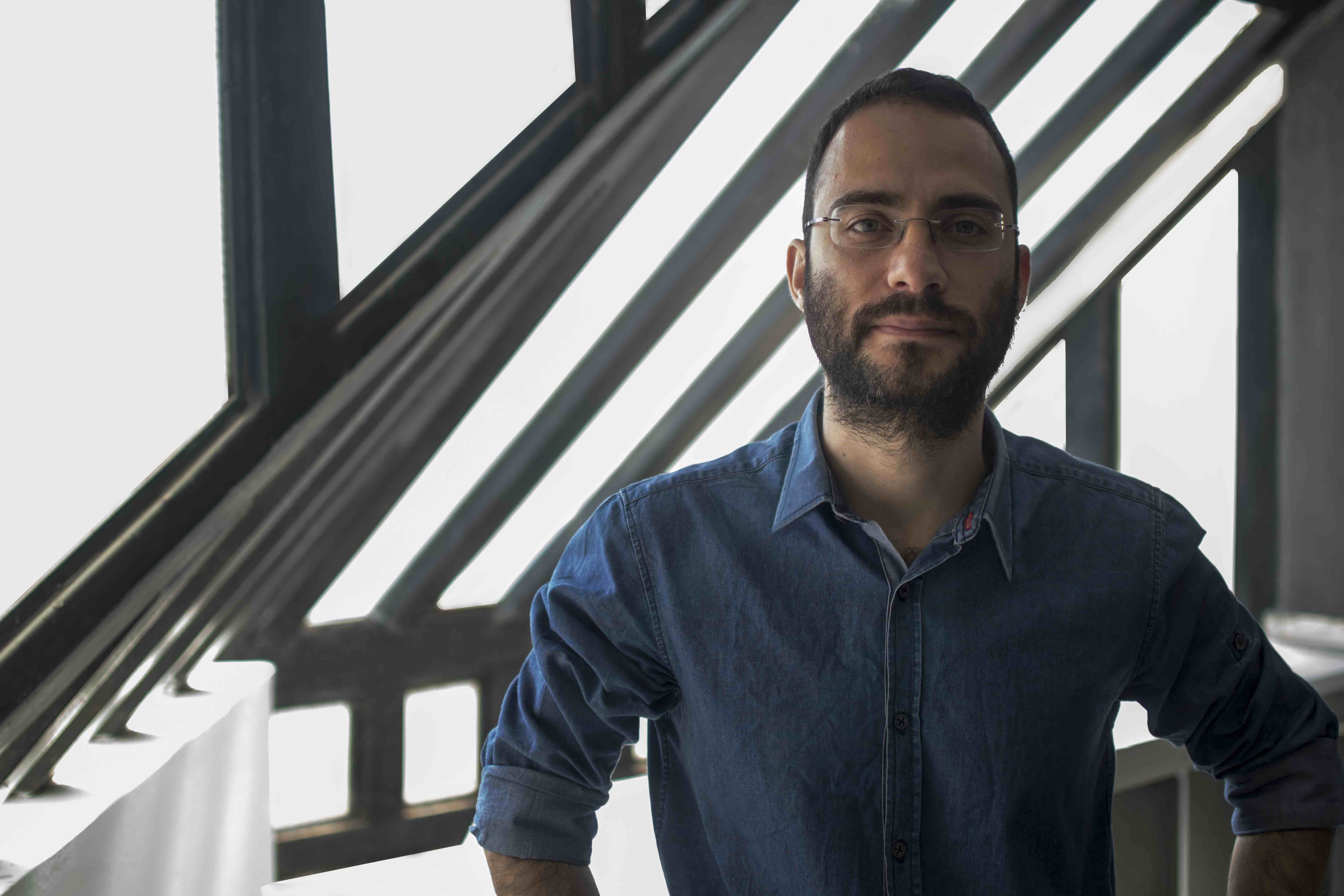}
\end{wrapfigure}
Nikolaos Bekiaris-Liberis received the Ph.D. degree in Aerospace Engineering from the University of California, San Diego, in 2013. From 2013 to 2014 he was a postdoctoral researcher at the University of California, Berkeley and from 2014 to 2017 he was a research associate and adjunct professor at Technical University of Crete, Greece. Dr. Bekiaris-Liberis is currently a Marie Sklodowska-Curie Fellow at the Dynamic Systems \& Simulation Laboratory, Technical University of Crete. He has coauthored the SIAM book {\em Nonlinear Control under Nonconstant Delays}. His interests are in delay systems, distributed parameter systems, nonlinear control, and their applications. 

Dr. Bekiaris-Liberis was a finalist for the student best paper award at the 2010 ASME Dynamic Systems and Control Conference and at the 2013 IEEE Conference on Decision and Control. He received the Chancellor's Dissertation Medal in Engineering from UC San Diego, in 2014. Dr. Bekiaris-Liberis is the recipient of a 2017 Marie Sklodowska-Curie Individual Fellowship Grant.\\[0.1cm]


\begin{wrapfigure}{L}{0.18\textwidth}
\vspace{-4mm}
\includegraphics[width=0.18\textwidth]{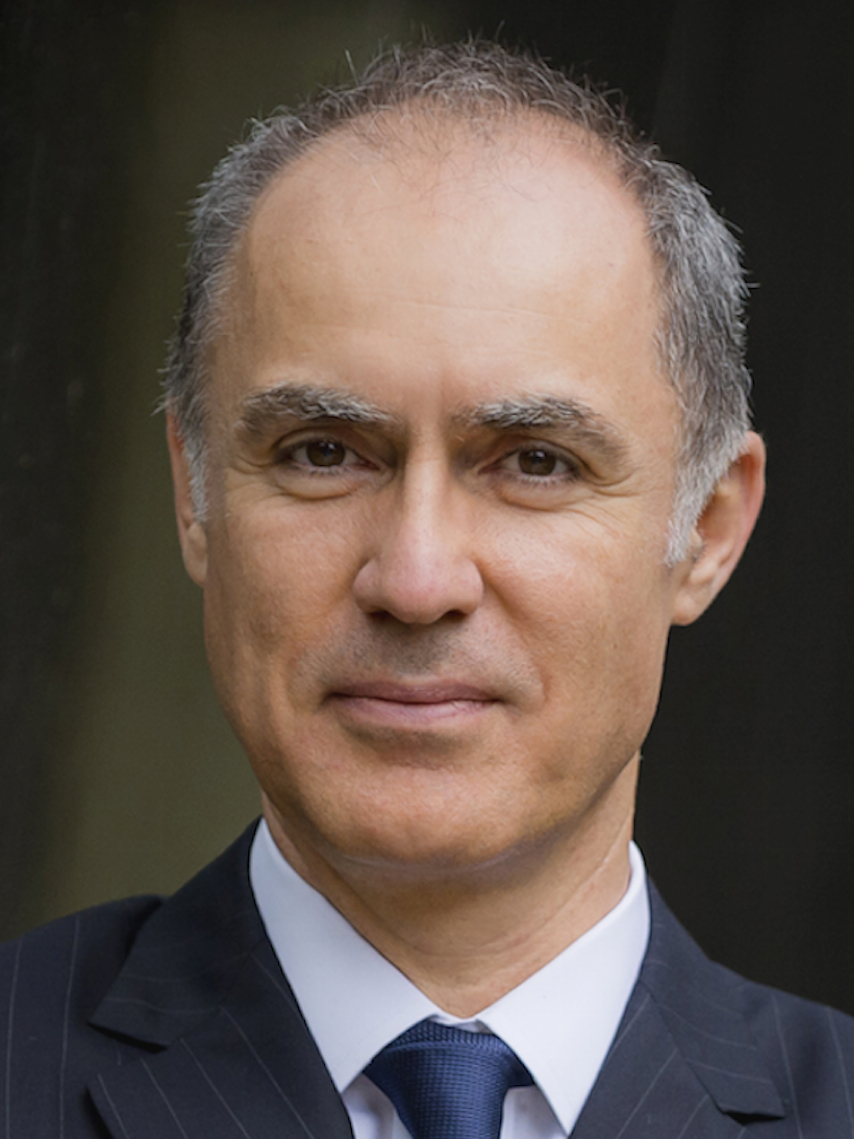}
\end{wrapfigure}
Miroslav Krstic is Distinguished Professor of Mechanical and Aerospace Engineering, holds the Alspach endowed chair, and is the founding director of the Cymer Center for Control Systems and Dynamics at UC San Diego. He also serves as Associate Vice Chancellor for Research at UCSD. As a graduate student, Krstic won the UC Santa Barbara best dissertation award and student best paper awards at CDC and ACC. Krstic is Fellow of IEEE, IFAC, ASME, SIAM, and IET (UK), Associate Fellow of AIAA, and foreign member of the Academy of Engineering of Serbia. He has received ASME Oldenburger Medal, ASME Nyquist Lecture Prize, ASME Paynter Outstanding Investigator Award, the PECASE, NSF Career, and ONR Young Investigator awards, the Axelby and Schuck paper prizes, the Chestnut textbook prize, and the first UCSD Research Award given to an engineer. Krstic has also been awarded the Springer Visiting Professorship at UC Berkeley, the Distinguished Visiting Fellowship of the Royal Academy of Engineering, the Invitation Fellowship of the Japan Society for the Promotion of Science, and honorary professorships from four universities in China. He serves as Senior Editor in IEEE Transactions on Automatic Control and Automatica, as editor of two Springer book series, and has served as Vice President for Technical Activities of the IEEE Control Systems Society and as chair of the IEEE CSS Fellow Committee. Krstic has coauthored twelve books on adaptive, nonlinear, and stochastic control, extremum seeking, control of PDE systems including turbulent flows, and control of delay systems.

}




\end{document}